\documentclass[twoside,parskip]{article}

	\usepackage[a4paper, left=4cm, right=4cm, top=4cm, bottom=4cm]{geometry}	
	\usepackage{amsmath}	
	\usepackage{amsthm}		
	\usepackage{amssymb}	
	\usepackage[all,cmtip]{xy}	
	\usepackage{amscd}		
	\usepackage{enumitem}	
	\usepackage{ifmtarg}	
	\usepackage{graphicx}	
	\usepackage{tocloft}	
	\usepackage[parfill]{parskip}	
	\usepackage{fancyhdr}		
	\usepackage{slashed}	
	\usepackage[
		sorting=nyt,             
		maxnames=5,              
		minnames=3,              
		backend=bibtex8,         
	]{biblatex}	
	\usepackage[pdftex,plainpages=false,pdfpagelabels]{hyperref}	
	\usepackage{color}	
	\usepackage{cleveref}	

	\bibliography{refs}
	
	\definecolor{linkcolor}{RGB}{00,10,138}
	\hypersetup{
		colorlinks=true,
		breaklinks=true,
		linkcolor=linkcolor,
		menucolor=linkcolor,
		citecolor=linkcolor,
		filecolor=linkcolor,
		urlcolor=linkcolor,
		frenchlinks=false
	}
	
	\begingroup
    \makeatletter
    \@for\theoremstyle:=definition,remark,plain\do{%
        \expandafter\g@addto@macro\csname th@\theoremstyle\endcsname{%
            \addtolength\thm@preskip\parskip
            }%
        }
	\endgroup
	
	\newcommand{\arsinh}{\operatorname{arsinh}}
    
	\newcommand{\ash}{\operatorname{a}} 
	\newcommand{\Conf}{\operatorname{Conf}}
	\newcommand{\Dirac}{\slashed{D}} 
	\newcommand{\Diff}{\operatorname{Diff}} 
	\newcommand{\dom}{\operatorname{dom}} 
	\newcommand{\ev}{\operatorname{ev}}	
	\newcommand{\GL}{\operatorname{GL}} 
	\newcommand{\GLp}{\GL^+} 
	\newcommand{\GLtp}{\widetilde{\GL}^+} 
	\newcommand{\image}{\operatorname{im}} 
	\newcommand{\id}{\operatorname{id}} 
	\newcommand{\ind}{\operatorname{index}} 
	\newcommand{\Mon}{\operatorname{Mon}} 
	\newcommand{\SO}{\operatorname{SO}} 
	\newcommand{\spc}{\mathfrak{s}}	
    \newcommand{\sa}{\operatorname{sa}}	
	\newcommand{\spec}{\operatorname{spec}} 
	\newcommand{\spp}{\mathfrak{sp}}	
	\newcommand{\specfl}{\operatorname{sf}} 
	\newcommand{\Spin}{\operatorname{Spin}} 
	\newcommand{\spin}{\operatorname{spin}} 
	\newcommand{\Rm}{\mathcal{R}} 
	
	\newcommand{\N}{\ensuremath{\mathbb{N}}\,}
	\newcommand{\Z}{\ensuremath{\mathbb{Z}}\,}
	\newcommand{\R}{\ensuremath{\mathbb{R}}\,}
	\newcommand{\C}{\ensuremath{\mathbb{C}}\,}
	
	\newcommand{\T}{\ensuremath{\mathbb{T}}\,}
	\renewcommand{\S}{S}
	
	\numberwithin{equation}{section}

	\newcounter{stepcounter}
	\makeatletter
	\newcommand{\stepitem}[1][]{ %
	\addtocounter{stepcounter}{1}
	\@ifmtarg{#1}{\item[{\sc{Step} \arabic{stepcounter}}:]}{\item[{\sc{Step} \arabic{stepcounter}} \normalfont{(#1)}:]}}
	\makeatother
	\newenvironment{steps}{\begin{description}[leftmargin=0ex,noitemsep,nosep]}{\end{description}\setcounter{stepcounter}{0}}
	
	\newcommand{\jeq}[2]{\stackrel{\text{#1}}{#2}}

	\newcommand{\DefMap}[4]{
	\begin{align*}
	\begin{array}{rcl}
	#1 & \to & #2 \\
	#3 & \mapsto & #4 
	\end{array} 
	\end{align*}
	}
	
	\theoremstyle{definition}
	\newtheorem{MainThm}{Main Theorem}
	\swapnumbers
	
	\newtheorem{Def}{Definition}[section]
	\newtheorem{Lem}[Def]{Lemma}
	\newtheorem{Thm}[Def]{Theorem}
	\newtheorem{Cor}[Def]{Corollary}
	\newtheorem{Rem}[Def]{Remark}

	\newtheorem{Not}[Def]{Notation}

	\newenvironment{Prf}{\begin{proof}[\normalfont \bfseries  Proof]}{\end{proof}}
	\newenvironment{Prf1}{\begin{proof}[\normalfont \bfseries  Proof of Main Theorem 1]}{\end{proof}}
	\newenvironment{Prf2}{\begin{proof}[\normalfont \bfseries  Proof of Main Theorem 2]}{\end{proof}}

\begin{document}

	\title{Continuity of Dirac Spectra
	}

	\author{Nikolai Nowaczyk
		\thanks{
			Universit\"at Regensburg,  
			Fakult\"at f\"ur Mathematik , 
			Universit\"atsstr. 31,  
			93040 Regensburg,
			Tel.: +49 941 9435692, 
			E-Mail: mail@nikno.de
		}
	}


	\date{2013-03-26}

	\maketitle

		\textbf{Abstract.} It is a well-known fact that on a bounded spectral interval the Dirac spectrum can be described locally by a non-decreasing sequence of continuous functions of the Riemannian metric. In the present article we extend this result to a global version. We think of the spectrum of a Dirac operator as a function $\Z \to \R$ and endow the space of all spectra with an $\arsinh$-uniform metric. We prove that the spectrum of the Dirac operator depends continuously on the Riemannian metric. As a corollary, we obtain the existence of a non-decreasing family of functions on the space of all Riemannian metrics, which represents the entire Dirac spectrum at any metric. We also show that in general these functions do not descend to the space of Riemannian metrics modulo spin diffeomorphisms due to spectral flow. \\ $ $ $ $ 


	\tableofcontents
	\thispagestyle{empty} 



\section{Introduction and Statement of the Results} \label{SctIntro}

For this entire article let $M^m$ be a smooth compact oriented spin manifold. Sticking to \cite[sec. 2]{BaerGaudMor} we fix a topological spin structure $\Theta:\GLtp M \to \GLp M$ on $M$. We denote by $\Rm(M)$ the space of all Riemannian metrics on $M$ endowed with $\mathcal{C}^1$-topology. For any metric $g \in \Rm(M)$ one obtains a metric spin structure $\Theta^g:\Spin^g M \to \SO^g M$ and an associated spinor bundle $\Sigma^g M$. The Dirac operator $\Dirac^g$ can be thought of as an unbounded operator $\Gamma_{L^2}(\Sigma^g M) \to \Gamma_{L^2}(\Sigma^gM)$ with domain $\Gamma_{H^1}(\Sigma^g M)$, where $H^1$ denotes the first order Sobolev space. For a more comprehensive introduction to spin geometry see \cite{LM}, \cite{FriedSpinGeom}.

Unfortunately one cannot directly compare the Dirac operators $\Dirac^g$ and $\Dirac^h$ for different metrics $g,h \in \Rm(M)$, because they are not defined on the same spaces. This problem has been discussed at length and solved in various other articles before, c.f. \cite{BaerGaudMor}, \cite{BourgGaud}, \cite{MaierGen}. The idea is to construct an isometry $\bar \beta_{g,h}:\Gamma_{L^2}(\Sigma^g M) \to \Gamma_{L^2}(\Sigma^h M)$ between the different spinor bundles and pull back the operator $\Dirac^h$ to an operator $\Dirac^h_g$ on the domain of $\Dirac^g$. This enables us to think of the Dirac operators as operators that depend continuously on the metric $g$ (the precise results are cited in \cref{ThmSpinorIdentification} later). It is therefore natural to ask, if and in what sense the spectrum of the Dirac operator also depends continuously on the metric. In the present article we will investigate this problem and present a solution. 

Every Dirac operator $\Dirac^g$ is a self-adjoint elliptic first order differential operator. It is a well-known fact (see for instance \cite[Thm. 5.8]{LM}) that the spectrum $\spec \Dirac^g$ is a subset of the real line that is closed, discrete and unbounded from both sides. The elements of $\spec \Dirac^g$ consist entirely of eigenvalues of finite multiplicity. Intuitively we would like to enumerate the eigenvalues from $- \infty$ to $+ \infty$ using $\Z$ as an index set by a non-decreasing sequence (in the entire paper we will always count eigenvalues with their geometric multiplicity). The problem is that this is not well-defined, because it is unclear which eigenvalue should be the ''first'' one. Formally we can avoid this problem as follows.

\begin{Def}
	For any $g \in \Rm(M)$ let $\spc^g:\Z \to \R$ be the unique non-decreasing function such that $\spc^{g}(\Z)=\spec \Dirac^{g}$, 
	\begin{align*}
		\forall \lambda \in \R: \dim \ker(\Dirac^g - \lambda) = \sharp (\spc^g)^{-1}(\lambda),
	\end{align*}
	 and $\spc^g(0)$ is the first eigenvalue $\geq 0$ of $\Dirac^g$. 
\end{Def}

Then $\spc^g$ is well-defined, but as it will turn out, the requirement that $\spc^g(0)$ is the first eigenvalue $\geq 0$ has some drawbacks. Namely, the map $g \mapsto \spc^g(j)$, $j \in \Z$, will not be continuous in general, see \cref{RemIntuitionDisc}. To obtain a more natural notion, we define the following.

\begin{Def}[$\Mon$ and $\Conf$]
	\label{DefMonConf}
	Define
	\begin{align*} 
		\Mon := \{ u:\Z \to \R \mid u \text{ is non-decreasing and proper} \} \subset  \R^\Z.
	\end{align*} 
	The group $(\Z,+)$ acts canonically on $\Mon$ via shifts, i.e.
	\begin{align} \label{EqDeftau}
	\begin{split}
		\begin{array}{rcl}
		\tau: \Mon \times \Z & \to & \Mon \\
		(u,z) & \mapsto & (j \mapsto (u.z)(j):=u(j+z))
		\end{array} 
	\end{split}
	\end{align}
	and the quotient
	\begin{align*}
		\Conf := \Mon / \Z
	\end{align*} 
	is called the \emph{configuration space}. Let $\pi:\Mon \to \Conf$, $u \mapsto \bar u$, be the quotient map.
\end{Def}

By construction $\spc^g \in \Mon$ and $\overline{\spc}^g := \pi (\spc^g) \in \Conf$. This defines maps
\begin{align*}
\spc:\Rm(M) \to \Mon, && \overline{\spc}:\Rm(M) \to \Conf.
\end{align*}
We would like to claim that $\overline{\spc}$ is continuous. To make formal sense of this, we introduce a topology on $\Mon$ and $\Conf$.

\begin{Def} [$\arsinh$-topology]
	\label{DefArsinhTopology}
	The topology induced by the metric $d_a$ defined by
	\begin{align*}
		\forall u,v \in \R^{\Z}: d_a(u,v) := \sup_{j \in \Z}{|\arsinh(u(j)) - \arsinh(v(j))}| \in [0,\infty]
	\end{align*}
	on $\R^\Z$ is called \emph{$\arsinh$-topology}. The group action $\tau$ acts by isometries with respect to $d_a$ and the quotient topology on $\Conf$ is induced by the metric $\bar d_a$ described by
	\begin{align} \label{EqDefbarda}
		\forall u \in \bar u, v \in \bar v \in \Conf: \bar d_a(\bar u, \bar v) = \inf_{j \in \Z}{d_a(u,v.j)}.
	\end{align}
	Using this metric on the quotient is common in metric geometry, c.f. \cite[Lemma 3.3.6]{Burago}.
\end{Def}

This allows us to formulate our main result.

\begin{MainThm}
	\label{MainThmSpec}
	The map $\overline{\spc} = \pi \circ \spc$ admits a lift $\widehat{\spc}$ against $\pi$ such that
	\begin{align}
        \label{EqMainThmSpc}
        \begin{split}
            \xymatrix{
                & (\Mon,d_a)
                    \ar[d]^-{\pi} 
                \\
                (\Rm(M),\mathcal{C}^1)
                    \ar[r]^-{\overline{\spc}}
                    \ar[ur]^-{\widehat{\spc}}
                & (\Conf, \bar d_a)}
        \end{split}
	\end{align}
	is a commutative diagram of topological spaces.
\end{MainThm}

From this one can quickly conclude the following claim, which is maybe a little more intuitive.

\begin{MainThm}
	\label{MainThmFun}
	There exists a family of functions $\{ \lambda_j \in \mathcal{C}^0(\Rm(M),\R) \}_{j \in \Z}$ such that for all $g \in \Rm(M)$ the sequence $(\lambda_j(g))_{j \in \Z}$ represents all the eigenvalues of $\Dirac^g$ (counted with multiplicities). In addition the sequence $\arsinh(\lambda_j)$ is equicontinuous and non-decreasing, i.e. all $g \in \Rm(M)$ satisfy $\lambda_j(g) \leq \lambda_k(g)$, if $j \leq k$.
\end{MainThm}

\begin{Prf2} Clearly the evaluation $\ev_j:(\Mon,d_a) \to \R$, $u \mapsto u(j)$, is a continuous map for any $j \in \Z$. Consequently, by \cref{MainThmSpec} the functions $\{ \lambda_j := \ev_j \circ \widehat{\spc} \}_{j \in \Z}$ are continuous and satisfy the assertion of \cref{MainThmFun}.
\end{Prf2}

\begin{figure}[t] 
	\begin{center}
		\includegraphics[scale=1, clip=true, trim=150 520 310 120]{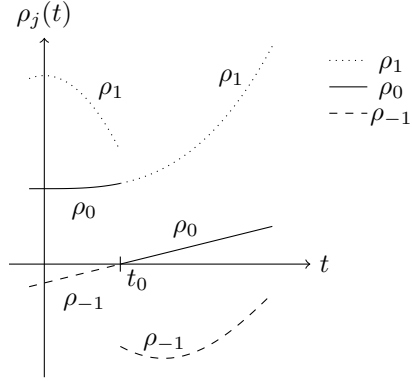}
		\caption{A null of $\rho_{-1}$ at $t_0$ can cause discontinuities at $t_0$ in all $\rho_j$.}
		\label{FigJumps}
	\end{center}
\end{figure}

\begin{Rem}[intuitive explanation]
	\label{RemIntuitionDisc}
	A choice of $\{\lambda_j\}_{j \in \Z}$ of functions representing all the Dirac eigenvalues that depends continuosly on $g \in \Rm(M)$ is a more subtle problem than one might think. The functions induced by $\spc$ (let's call these $\rho_j := \ev_j \circ \spc$ for the moment) are not continuous in general. To see this imagine a continuous path of metrics $(g_t)_{t \in \R}$ and consider $\rho_j:\R \to \R$ as functions of $t$, see \Cref{FigJumps}. Since $\rho_0(t)$ is the first eigenvalue $\geq 0$ of $\Dirac^{g_t}$, this function will have a jump at point $t_0$ where $\rho_0(t_0)>1$ and $\rho_{-1}(t_0)=0$. This can cause discontinuities in all the other functions $\rho_j$ as well. \\
	However for any $k \in \Z$ the sequence $\rho_j' := \rho_{j+k}$, $j \in \Z$, gives another enumeration of the spectrum. Intuitively \cref{MainThmSpec} states that if one uses this freedom in the enumeration of the eigenvalues at each metric in the ''right'' way, one obtains a globally well-defined family of continuous functions representing all the Dirac eigenvalues.
\end{Rem}

The rest of this paper is organized as follows: After a short review of some fundamental results in \Cref{SctFundamental} and a slight generalization of our notation in \Cref{SctNotSpect} the main part of the paper will be \Cref{SctProofMainThm}, which is devoted to build up technical results that are needed for the proof of \cref{MainThmSpec}. Finally in \Cref{SecModuliSpecFl} we will investigate to what extent the functions $\widehat{\spc}$ and $\overline{\spc}$ descend to certain quotients of $\Rm(M)$ called \emph{moduli spaces}. Our central result will be that there exists an obstruction, the \emph{spectral flow}, for $\widehat{\spc}$ to descend onto $\Rm(M) / \Diff^{\spin}(M)$. This will be made precise in \cref{DefSpinDiffeo} and \cref{LemSpecFlowDef} and the cetral result will be stated in \cref{MainThmFlow}.

Using these results the actual proof of \cref{MainThmSpec} becomes very short.

\begin{Prf1} By \cref{ThmSpecContArsinh} the map $\bar \spc:(\Rm(M),\mathcal{C}^1) \to (\Conf, \bar d_a)$ is continuous. By \cref{ThmPiCoveringArsinh} the map $\pi:(\Mon,d_a) \to (\Conf,\bar d_a)$ is a covering map. Since $\Rm(M)$ is path-connected, locally path-connected and simply connected, the claim follows from the Lifting Theorem of Algebraic Topology.
\end{Prf1}

\begin{Rem}[uniqueness]
	From this proof we conclude that the lift $\widehat{\spc}$ is not unique. In fact there are $\Z$ possibilities of how to lift $\overline{\spc}$ against $\pi$. One can use this freedom to arrange that $\widehat{\spc}^{g_0} = \spc^{g_0}$ for one fixed $g_0 \in \Rm(M)$.
\end{Rem}

\section{Fundamental Results} \label{SctFundamental}

\begin{Thm}[{identification of spinor bundles, c.f. \cite{MaierGen}}]
	\label{ThmSpinorIdentification}
	Let $g \in \Rm(M)$ be a fixed metric. For every $h \in \Rm(M)$ there exists an isometry of Hilbert spaces $\bar \beta_g^h: \Gamma_{L^2}(\Sigma^gM) \to \Gamma_{L^2}(\Sigma^hM)$ such that the operator
	\begin{align*}
	\Dirac^h_g := \bar \beta^g_h \circ \Dirac^h \circ \bar \beta_g^h: \Gamma_{L^2}(\Sigma^gM) \to \Gamma_{L^2}(\Sigma^gM)
	\end{align*}
	is closed, densely defined on $\Gamma_{H^1}(\Sigma^g M)$, isospectral to $\Dirac^h$, and the map
	\begin{align*}
	\Dirac_g:\Rm(M) \to B(\Gamma_{H^1}(\Sigma^gM),\Gamma_{L^2}(\Sigma^gM)), && h \mapsto \Dirac_g^h,
	\end{align*}
	is continuous. (Here $B(\_)$ denotes the space of bounded linear operators endowed with the operator norm.)
\end{Thm}

The following theorem is formally not needed for the proof of \cref{MainThmSpec}. It is nevertheless worth mentioning to give an impression of what is already well-known about the continuity of Dirac spectra. It implies that a bounded spectral interval of the Dirac operator can be described locally by continuous functions. Consequently \cref{MainThmSpec} can be thought of as a global analogue of this local result.

\begin{Thm}[{\cite[Prop. 7.1]{BaerMetrHarmonSpin}}]
	\label{ThmBaerSpecC1Conti}
	Let $(M,g)$ be a closed Riemannian spin manifold with Dirac operator $\Dirac^g$ having spectrum $\spec \Dirac^g$. Let $\Lambda > 0$ such that $-\Lambda,\Lambda \notin \spec \Dirac^g$ and enumerate
	\begin{align*}
		\spec \Dirac^g \cap \mathopen] -\Lambda,\Lambda \mathclose[ = \{ \lambda_1 \leq \lambda_2 \leq \ldots \leq \lambda_n \}.
	\end{align*}
	For any $\varepsilon > 0$ there exists a $\mathcal{C}^1$-neighbourhood $U$ of $g$ such that for any $g' \in U$ 
	\begin{enumerate}
		\item $\spec \Dirac^{g'} \cap \mathopen] -\Lambda, \Lambda \mathclose[ = \{ \lambda'_1 \leq \ldots \leq \lambda'_n \}$,
		\item $\forall 1 \leq i \leq n: |\lambda_i - \lambda_i'| < \varepsilon$.
	\end{enumerate}
\end{Thm}

\section{Families of Discrete Operators} \label{SctNotSpect}
For the proof of \cref{MainThmSpec} we need the notion of the function $\spc^g$ for operators that are slightly more general than Dirac operators. In this section we introduce the necessary definitions and notation conventions. Let $X,Y$ be complex Banach spaces. We denote by $C(X,Y)$ the space of unbounded operators $T:X \supset \dom(T) \to Y$. Let $B(X,Y)$ denote the bounded operators $X \to Y$. We set $C(X):=C(X,X)$ and $B(X):=B(X,X)$. The spectrum of $T$ is denoted by $\spec T \subset \C$.

\begin{Def}[discrete operator]
	\label{DefDiscreteOperator}
	An operator $T \in C(X)$ is \emph{discrete}, if $\spec T \subset \R$ is a closed discrete subset that is unbounded from both sides and consists solely of eigenvalues that are of finite multiplicity.
\end{Def}

\begin{Def}[ordered spectral function] 
	Let $T \in C(X)$ be discrete. The sequence $\spc_T \in \R^\Z$ uniquely defined by the properties
	\begin{enumerate}
		\item $\spc_T(0) = \min \{\lambda \in \spec T \mid \lambda \geq 0 \}$.
		\item $\forall i,j \in \Z: i \leq j \Longrightarrow \spc_T(i) \leq \spc_T(j)$.
		\item $\forall \lambda \in \R: \sharp (\spc_T)^{-1}(\lambda)= \dim \ker (T - \lambda) $.
	\end{enumerate}
	is the \emph{ordered spectral function of $T$}. 
\end{Def}

\begin{Def}[spectral parts] 
	Let $T \in C(X)$ be discrete. To denote parts of the ordered spectrum, we introduce the following notation: If $I \subset \R$ is an interval, then $(\spc_T)^{-1}(I) = \{k, k+1, \ldots, l\}$ for some $k, l \in \Z$, $k \leq l$. The sequence
	\begin{align*} 
		\spp_T(I) := (\spc_T(i))_{k \leq i \leq l},
	\end{align*}
	is the \emph{spectral part of $T$ in $I$}. 
\end{Def}

\begin{Def}[discrete family]
	Let $E$ be any set. A map $T:E \to C(X)$ is a \emph{discrete family}, if for any $e \in E$ the operator $T_e$ is discrete in the sense of \cref{DefDiscreteOperator}. We obtain a function
	\DefMap{\spc_T:E}{\R^\Z}{e}{\spc_T^e:=\spc_{T(e)}.}
	Analogously, we set $\spp_T^e := \spp_{T(e)}$.
\end{Def}

\begin{Rem}[family of Dirac operators]
	In view of \cref{ThmSpinorIdentification} we can apply the above in particular to Dirac operators. Namely, we fix any $g \in \Rm(M)$ and set $X:=L^2(\Sigma^g M)$ and $E:=\Rm(M)$. Then $h \mapsto \Dirac^h_g$ is a discrete family. We will suppress its name in notation and just write $\spc^h_g$ for the ordered spectral function of $\Dirac^h_g$. Since $\Dirac^h$ and $\Dirac^h_g$ are isospectral, we can ignore the reference metric entirely and just write $\spc^h$.
\end{Rem}

\section{Proof of Main Theorem 2} \label{SctProofMainThm}

In this section we carry out the details of the proof of \cref{MainThmSpec}. The idea to construct the $\arsinh$-topology in the first place is inspired by a paper from John Lott, c.f. \cite[Theorem 2]{Lott}. The arguments require some basic notions from analytic pertubation theory. We cite and slightly modify some results by Kato, c.f. \cite{kato}. Applying analytic pertubation theory to families of Dirac operators is a technique that is also used in other contexts, c.f. \cite{BaerGaudMor}, \cite{BourgGaud}, \cite{AndreasDiss}. \\
Let $X$,$Y$ be complex Banach spaces and let $X'$ be the topological dual space of $X$. For any operator $T$ we denote its adjoint by $T^*$. Let $\Omega \subset \C$ be an open and connected subset. Recall that a function $f:\Omega \to X$ is \emph{holomorphic}, if for all $\zeta_0 \in \Omega$
\begin{align*}
	f'(\zeta_0) := \lim_{\zeta \to \zeta_0}{\tfrac{f(\zeta)-f(\zeta_0)}{\zeta-\zeta_0}}
\end{align*}
exists in $(X,\|\_\|_X)$. A family of operators $T:\Omega \to B(X,Y)$ is \emph{bounded holomorphic}, if $T$ is a holomorphic map in the sense above. To treat the unbounded case, the following notions are crucial.

\begin{Def}[holomorphic family of type (A)]
	A family of operators $T:\Omega \to C(X,Y)$, $\zeta \mapsto T_\zeta$, is \emph{holomorphic of type (A)}, if the domain $\dom(T_\zeta) =: \dom(T)$ is independent of $\zeta$ and for any $x \in X$ the map $\Omega \to Y$, $\zeta \mapsto T_\zeta x$, is holomorphic.
\end{Def}

\begin{Def}[self-adjoint holomorphic family of type(A)]
	A family $T:\Omega \to C(H)$ is \emph{self-adjoint holomorphic of type (A)}, if it is holomorphic of type (A), $H$ is a Hilbert space, $\Omega$ is symmetric with respect to complex conjugation and
	\begin{align*}
		\forall \zeta \in \Omega: T^*_{\zeta} = T_{\bar \zeta}.
	\end{align*}
\end{Def}

These families are particularly important for our purposes due to the following useful theorem.

\begin{Thm}[{\cite[VII.\S 3.5, Thm. 3.9]{kato}}]
	\label{ThmAnalyticEVKato}
	Let $T:\Omega \to C(H)$ be a self-adjoint holomorphic family of type (A) and let $I \subset \Omega \cap \R$ be an interval. Assume that $T$ has compact resolvent. Then there exists a family of functions $\{\lambda_n \in \mathcal{C}^\omega(I,\R)\}_{n \in \N}$ and a family functions $\{u_n \in \mathcal{C}^\omega(I,H)\}$ such that for all $t \in I$, the $(\lambda_n(t))_{n \in \N}$ represent all the eigenvalues of $T_t$ counted with multiplicity, $T_tu_n(t)=\lambda_n(t)u_n(t)$, and the $(u_n(t))_{n \in \N}$ form a complete orthonormal system of $H$.
\end{Thm}

Derivatives of holomorphic families can be estimated using the following theorem.

\begin{Thm}[{\cite[VII.\S 2.1, p.375f]{kato}}]
	\label{ThmHolAFamilyDerGrowth}
	Let $T:\Omega \to C(X,Y)$ be a holomorphic family of type (A). For any $\zeta \in \Omega$ define the operator
	\begin{align*}
		T'_{\zeta}: \dom(T) \to Y, && u \mapsto T'_{\zeta}u := \tfrac{d}{d\zeta}(T_{\zeta}u).
	\end{align*}
	Then $T'$ is a map from $\Omega$ to the unbounded operators $X \to Y$ (but $T'_{\zeta}$ is in general not closed). For any compact $K \subset \Omega$ 
	there exists $C_K>0$ such that 
	\begin{align*}
		\forall \zeta \in K: \; \forall u \in \dom(T): \; \|T'_{\zeta}u\|_Y \leq C_K(\|u\|_X + \|T_{\zeta}u\|_Y).
	\end{align*}
	If $\zeta_0 \in K$ is arbitrary, $Z:=\dom(T)$ and $\|u\|_Z := \|u\|_X + \|T_{\zeta_0}u\|_Y$, then $C_K := \alpha_K^{-1} \beta_K$ does the job, where
	\begin{align} \label{EqHolAFamilyDerGrowthConstantSpec}
		\alpha_K := \inf_{\zeta \in K}{ \inf_{\|u\|_Z=1}{\|u\|_X + \|T_{\zeta}u\|_Y}}, &&
		\beta_K := \sup_{\zeta \in K}{\|T'_{\zeta}\|_{B(Z,Y)}}.
	\end{align}
\end{Thm}

This can be used to prove the following about the growth of eigenvalues. 

\begin{Thm}[{\cite[VII.\S 3.4, Thm. 3.6]{kato}}]
	\label{ThmHolAFamilyEVGrowth}
	Let $T:\Omega \to C(H)$ be a self-adjoint holomorphic family of type (A).  Let $I \subset \Omega \cap \R$ be a compact interval and let $J \subset I$ be open. Assume that $\lambda \in \mathcal{C}^\omega(J,\R)$ is an \emph{eigenvalue function}, i.e. for all $ t \in J$ the value $\lambda(t)$ is an eigenvalue of $T_t$. Then
	\begin{align} \label{EqHolAFamilyConstantExp}
		\forall t,t_0 \in J: \; |\lambda(t) - \lambda(t_0)| \leq (1 + |\lambda(t_0)|)(\exp(C_I |t-t_0|) - 1),
	\end{align}
	where $C_I$ is the constant from \cref{ThmHolAFamilyDerGrowth}. 
\end{Thm}

The important consequence is that \eqref{EqHolAFamilyConstantExp} can be reformulated in terms of the $\arsinh$-topology. This is the first of the three main parts of the technical work towards \cref{ThmSpecContArsinh}.

\begin{Cor}[Growth of eigenvalues]
	\label{CorHolAFamilyEVGrowth}
	In the situation of \cref{ThmHolAFamilyEVGrowth} the following holds in addition: For any $t_0 \in I$ and $\varepsilon > 0$ there exists $\delta > 0$ such that for all $t \in I_{\delta}(t_0) \cap J$ and all eigenvalue functions $\lambda \in \mathcal{C}^\omega(J,\R)$
	\begin{align} \label{EqHolAFamilyDerGrowthArsinh}
		|\arsinh(\lambda(t)) - \arsinh(\lambda(t_0))| < \varepsilon.
	\end{align}
	There exist universal constants (i.e. independent of the family $T$) $C_1,C_2 > 0$ such that 
	\begin{align} \label{EqHolAFamilyDerGrowthDelta}
		\delta := C_I^{-1} \ln(\min(C_1, \varepsilon C_2)  + 1 )
	\end{align} 
	does the job.
\end{Cor}

\begin{Prf} $ $

	\begin{steps}
	\stepitem
		The function $\alpha:\R \to \R$, $t \mapsto \exp(C_I |t-t_0|) - 1$, is continuous and satisfies $\alpha(t_0)=0$. Notice that for $b > 0$ 
		\begin{align}
			\label{EqHolAFamilyEVGrowthDeltaCalc}
			|\alpha(t)| < b \Longleftrightarrow |t-t_0| < C_I^{-1} \ln(b+1).
		\end{align}
		In particular there exists $\delta_1 > 0$ such that 
		\begin{align} \label{EqHolAFamilyConstantExpDeltaQuarter}
		\forall t \in I_{\delta_1}(t_0): |\alpha(t)| < \tfrac{1}{4}
		\end{align}
		So let $t \in I_{\delta_1}(t_0)$. 
		
	\stepitem Setting $\lambda_0:=\lambda(t_0)$ we can reformulate \eqref{EqHolAFamilyConstantExp} by
	\begin{align} \label{HolAFamilyEVGrowthReform}
	\lambda_0 - (1 + |\lambda_0|)\alpha(t) < \lambda(t) < \lambda_0 +(1+|\lambda_0|) \alpha(t).
	\end{align}
	Since
	\begin{align*}
		\lim_{|R| \to \infty}{\tfrac{|R|}{1+|R|}} = 1 ,
	\end{align*}
	and the convergence is monotonously increasing, there exists $R > 0$ such that
	\begin{align} \label{EqEVfunArsinhExQuotient}
		\forall |\eta| \geq R: \tfrac{1}{2} < \tfrac{|\eta|}{1 + |\eta|}.
	\end{align}
	Now assume $|\lambda_0| \geq R$. In case $\lambda_0 \geq R > 0$, we calculate
	\begin{align}
		\label{EqHolAFamilyEVGrowthBoundBelowPositive}
		\tfrac{\lambda_0}{1 + \lambda_0} > \tfrac{1}{2} \jeq{\eqref{EqHolAFamilyConstantExpDeltaQuarter}}{\geq} 2 \alpha(t)		
		\Longrightarrow \tfrac{1}{2} \lambda_0 \geq \alpha(t)(1+\lambda_0)
		\Longrightarrow \lambda_0 - \alpha(t)(1 +|\lambda_0|) \geq  \tfrac{1}{2} \lambda_0 .
	\end{align}
	Analogously, if  $\lambda_0 \leq -R < 0$, we calculate
	\begin{align}
		\label{EqHolAFamilyEVGrowthBoundBelowNegative}
		\lambda_0 + \alpha(t) (1 + |\lambda_0|) < \tfrac{1}{2} \lambda_0.
	\end{align}

	\stepitem
		Define the constants
		\begin{align*}
			C_0 := \sup_{t \in \R}{\tfrac{1 + |t|}{\sqrt{1 + t^2}}}, &&
			C_1 := \tfrac{1}{4}, &&
			C_2 := \min \left( \tfrac{1}{R+1}, \tfrac{1}{2C_0} \right)
		\end{align*}
		and set 
		\begin{align*}
		\delta_2 := C_I^{-1} \ln(\min(C_1, \varepsilon C_2)  + 1 ).
		\end{align*}
		By \eqref{EqHolAFamilyEVGrowthDeltaCalc} this implies 
		\begin{align} \label{EqHolAFamilyEVGrowthDeltaInp}
			\forall t \in I_{\delta_2}(t_0): \alpha(t) < \min(C_1,\varepsilon C_2) \leq \varepsilon C_2.
		\end{align}
		So let $t \in I_{\delta_2}(t_0)$ be arbitrary and set $c_\pm := \lambda_0 \pm (1 + |\lambda_0|)\alpha(t) $. It follows from the Taylor series expansion of $\arsinh$ that there exists $\xi \in [\lambda_0,c_+]$ such that
		\begin{align} \label{EqHolAFamilyEVGrowthFinal} 
			\arsinh(c_+) - \arsinh(\lambda_0) 
			= \arsinh'(\xi) (1 + |\lambda_0|)\alpha(t) 
			=\frac{(1 + |\lambda_0|)}{\sqrt{1 + \xi^2}} \alpha(t).
		\end{align}
		Now in case $|\lambda_0| \leq R$, we continue this estimate by
		\begin{align*}
			\eqref{EqHolAFamilyEVGrowthFinal}
			\leq (1 + |\lambda_0|) \alpha(t)
			\leq (1 + R) \alpha(t) 
			\jeq{\eqref{EqHolAFamilyEVGrowthDeltaInp}}{<} \varepsilon
		\end{align*}
		In case $\lambda_0 \geq R$, we continue this estimate by
		\begin{align*}
			\eqref{EqHolAFamilyEVGrowthFinal}
			\leq \frac{(1 + |\lambda_0|)}{\sqrt{1 + \lambda_0^2}} \alpha(t)
			\leq C_0 \alpha(t)  
			\jeq{\eqref{EqHolAFamilyEVGrowthDeltaInp}}{<} \varepsilon.
		\end{align*}
		In case $\lambda_0 \leq -R$, we continue this estimate by
		\begin{align*}
			\eqref{EqHolAFamilyEVGrowthFinal}
			\leq \frac{(1 + |\lambda_0|)}{\sqrt{1 + c_+^2}} \alpha(t)
			\jeq{\eqref{EqHolAFamilyEVGrowthBoundBelowNegative}}{\leq} \frac{(1 + |\lambda_0|)}{\sqrt{1 + \tfrac{1}{4} \lambda_0^2}} \alpha(t)
			\leq 2 C_0 \alpha(t)  
			\jeq{\eqref{EqHolAFamilyEVGrowthDeltaInp}}{<} \varepsilon.
		\end{align*}
		Consequently, since $\arsinh$ is strictly increasing, in all cases we obtain 
		\begin{align*}
			\arsinh(\lambda(t))
			\jeq{\eqref{HolAFamilyEVGrowthReform}}{<}\arsinh(\lambda_0 + (1 + |\lambda_0|)\alpha(t))
			< \arsinh(\lambda_0) + \varepsilon.
		\end{align*}
		By a completely analogous argument, we obtain 
		\begin{align*}
			\arsinh(\lambda(t))
			> \arsinh(\lambda_0 - (1+|\lambda_0|)\alpha(t))
			\geq \arsinh(\lambda_0) - \varepsilon .
		\end{align*}
	\end{steps}
	This proves the claim.	
\end{Prf}

Now, the second step is to apply the preceeding result to discrete families.

\begin{Not} \label{NotTechnical}
	Since the following proof is a little technical, we abbreviate
	\begin{align*}
		\boxed{\ash := \arsinh}
	\end{align*}
	and set $d_a(x,y):=|\ash(x)-\ash(y)|$ for $x,y \in \R$. For any $\varepsilon > 0$ the $\varepsilon$-neighbourhoods of $x \in \R$ and $\varepsilon$-hulls of a set $S \subset \R$ will be denoted by 
	\begin{align*}
		I_\varepsilon(x) := \{ t \in \R \mid |t-x| < \varepsilon \}, &&
		I_\varepsilon(S) := \bigcup_{x \in S}{I_\varepsilon(x)}, &&
		I_\varepsilon^a(S) := I_\varepsilon(\ash(S)).
	\end{align*}
	For any $k,l \in \Z$, $k \leq l$, we define $[k,l]_\Z := [k,l] \cap \Z$.
\end{Not}

\begin{Cor}[Spectral Growth] \label{CorAnalyticGlobalGrowth}
	Let $\Omega \subset \C$ be open, $I \subset \Omega \cap \R$ be an interval, $T:\Omega \to C(H)$ be a discrete self-adjoint holomorphic family of type (A). For any $t_0 \in I$ and $\varepsilon > 0$ there exists $\delta > 0$ such that 
	\begin{align} \label{EqAnalyticGlobalGrowth}
		\forall t \in I_\delta(t_0) \cap I: \exists k \in \Z: \forall j \in \Z: d_a(\spc^{t_0}_T(j),\spc^{t}_T(j+k)) < \varepsilon.
	\end{align}
\end{Cor}

\begin{Prf} $ $ 
\begin{steps}
	\stepitem
		Certainly the family of eigenfunctions from \cref{ThmAnalyticEVKato} can be $\Z$-reindexed to a family $\{ \lambda_j \in \mathcal{C}^\omega(I,\R) \}_{j \in \Z}$ satisfying $\lambda_j(t_0)=\spc_T^{t_0}(j)$, $j \in \Z$. By \cref{CorHolAFamilyEVGrowth} 
		\begin{align*}
			\exists \delta > 0: \forall t \in I_\delta(t_0) \cap I: \forall j \in \Z: |\ash(\lambda_j(t_0)) - \ash(\lambda_j(t))| < \varepsilon .
		\end{align*}
		Fix any $t \in I_\delta(t_0) \cap I$ and let $\sigma:\Z \to \Z$ be the bijection satisfying $\spc_T^t(\sigma(j)) = \lambda_j(t)$. This implies 
		\begin{align} \label{EqAnalyticGlobalGrowthBij}
			\forall j \in \Z: d_a(\spc_T^{t_0}(j),\spc_T^t(\sigma(j)) < \varepsilon,
		\end{align}
		which ist almost \eqref{EqAnalyticGlobalGrowth}, except that $\sigma$ might not be given by a translation. In the next steps we will show that we may replace the bijection $\sigma$ by an increasing bijection $\tau$ that still satisfies \eqref{EqAnalyticGlobalGrowthBij}. Since every increasing bijection $\Z \to \Z$ is given by a translation $\tau^{(k)}:\Z \to \Z$, $z \mapsto z + k$, for some $k \in \Z$, this implies the claim.

	\stepitem
		For any $n \in \N$ consider the function $\sigma_n:=\sigma|_{[-n,n]_\Z}:[-n,n]_\Z \to \Z$. This function is injective and satisfies \eqref{EqAnalyticGlobalGrowthBij} for all $-n \leq j \leq n$. Furthermore setting
		\begin{align*}
			\spp^T_{t_0}([\lambda_{-n}(t_0), \lambda_n(t_0)]) &=: (\lambda_{-n}, \ldots, \lambda_{n})
		\end{align*}
		we obtain numbers $n',m' \in \Z$ such that the eigenvalues $\mu_j := \spc_T^{t}(j)$ satisfy
		\begin{align} \label{AnalyticLocalGrowthEVNumber}
			I_\varepsilon^a(\spp^T_{t}([\lambda_{-n},\lambda_{n}])) & = (\ash(\mu_{n'}), \ldots, \ash(\mu_{m'}) ), 
		\end{align}
		and we have the estimate
		\begin{align} \label{EqAnalyticGlobalGrowthLocalized} 
			\forall -n \leq j \leq n: & |\ash(\lambda_j) - \ash(\mu_{\sigma_n(j)}) | < \varepsilon.
		\end{align}
		We will show that $\sigma_n$ can be modified to an increasing injection $\tilde \sigma_n$ which satisfies $\image(\sigma_n) = \image(\tilde{\sigma}_n)$ and \eqref{EqAnalyticGlobalGrowthLocalized}. To that end choose any $-n \leq i < j \leq n$ and assume that $\sigma_n(j)<\sigma_n(i)$. Notice that by construction
		\begin{align}
			\label{EqAnalyticLocalGrowthMonotoneLambdaMu}
			i < j \Longrightarrow  \lambda_i \leq \lambda_j , && \sigma_n(j) < \sigma_n(i) \Longrightarrow \mu_{\sigma_n(j)} \leq \mu_{\sigma_n(i)}.
		\end{align}
		Define the function $\tilde \sigma_n$ by setting 
		\begin{align*}
			\tilde \sigma_n|_{\{-n, \ldots, n\} \setminus \{i,j\}} := \sigma_n, &&
		\tilde \sigma _n(i)=\sigma_n(j), &&
		\tilde \sigma_n(j)=\sigma_n(i).
		\end{align*}
		It is clear that $\tilde \sigma_n$ is still injective and $\image(\tilde \sigma_n) = \image(\sigma_n)$. To show that it still satisfies \eqref{EqAnalyticGlobalGrowthLocalized} we distinguish two cases, see figure \Cref{FigCones}. First consider the case that $\lambda_i = \lambda_j$. Then it follows automatically that
		\begin{align*}
			|\ash(\lambda_i) - \ash(\mu_{\tilde \sigma_n(i)})| = |\ash(\lambda_j) - \ash(\mu_{\sigma_n(j)})| < \varepsilon
		\end{align*}
		and the same for $j$. In case $\lambda_i \neq \lambda_j$ it follows that $\lambda_i < \lambda_j$. This implies
		\begin{align*}
			\ash(\lambda_i) - \varepsilon < \ash(\lambda_j) - \varepsilon < \ash(\mu_{\sigma_n(j)}) \jeq{\eqref{EqAnalyticLocalGrowthMonotoneLambdaMu}}{\leq} \ash(\mu_{\sigma_n(i)}) < \ash(\lambda_i) + \varepsilon < \ash(\lambda_j) + \varepsilon,
		\end{align*}
		hence
		\begin{align*}
			\mu_{\sigma_n(i)},\mu_{\sigma_n(j)} \in I_\varepsilon^a(\lambda_i) \cap I_\varepsilon^a(\lambda_j).
		\end{align*}
		In particular, this intersection is not empty. Consequently $\tilde \sigma_n$ satisfies \eqref{EqAnalyticGlobalGrowthLocalized}. By repeating this procedure for all index pairs $(i,j)$, $-n \leq i \leq n$, $i<j \leq n$, it follows that $\sigma_n$ can be modified finitely many times in this manner to obtain an increasing injection having the same image that still satisfies \eqref{EqAnalyticGlobalGrowthLocalized}. For simplicity denote this function still by $\tilde \sigma_n$ and define
		\DefMap{\tilde \tau_n:\Z}{\Z}{j}{\begin{cases}
		\tilde \sigma_n(j), & -n \leq j \leq n, \\
		\sigma(j), & \text{otherwise}.
		\end{cases}}
		This function is still bijective, still satisfies \eqref{EqAnalyticGlobalGrowthBij} and is increasing on $[-n,n]_\Z=:I_n$. Define $J_n:=\tilde{\tau}_n(I_n)$.

\begin{figure}[t] 
	\begin{center}
		\includegraphics[scale=0.9, clip=true, trim=90 680 230 70]{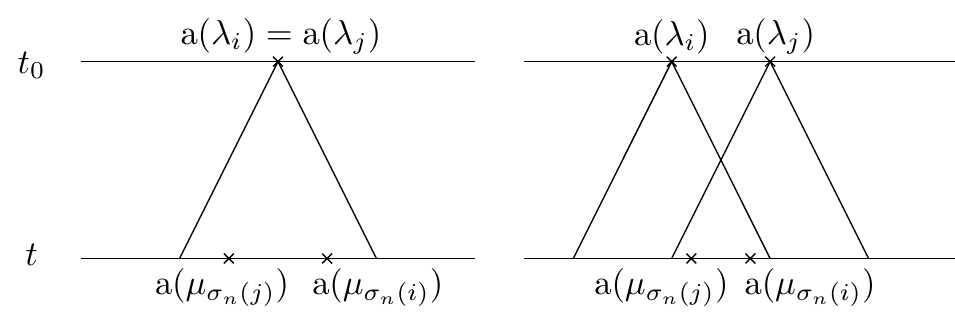}
		\caption{The two Possibilities for $\lambda_j$.}
		\label{FigCones}
	\end{center}
\end{figure}

	\stepitem 
		Unfortunately it might happen that $\tilde{\tau}_{n+1}|_{I_n} \neq \tilde{\tau}_n$. But due to \eqref{EqAnalyticGlobalGrowthLocalized} there exists $n_1$ such that all $n \geq n_1$ satisfy $\tilde{\tau}_n(I_1) \subset J_{n_1}$. Since there are only finitely many functions $I_{1} \to J_{n_1}$, there has to be at least one such function that occurs infinitely often in the sequence $\{\tilde{\tau}_n|_{I_1}\}_{n \in \N}$. So there exists an infinite subset $\N_1 \subset \N_0:=\N$ such that $\tilde{\tau}_n|_{I_1}$ is the same for all $n \in \N_1$. \\
		Now, the same holds for $I_2$: There exists $n_2 \geq n_1$ such that all $n \geq n_2$ satisfy $\tau_n(I_2) \subset J_{n_2}$. Again since there are only finitely many functions $I_2 \to J_{n_2}$, one of them has to occur infinitely often in the sequence $\{\tau_n|_{I_2}\}_{n \in \N_1}$. Consequently there exists an infinite subset $\N_2 \subset \N_1$ such that $\tau_n|_{I_2}$ is the same for all $n \in \N_2$. This process can be continued indefinitely for all the intervals $I_\nu$, $\nu \in \N$. Finally the function
		\DefMap{\tau:\Z}{\Z}{j}{\tilde{\tau}_n(j), \; j \in I_\nu, n \in \N_{\nu}}
		does the job: It is well-defined, satisfies \eqref{EqAnalyticGlobalGrowthBij}, it remains injective and it is surjective: Since the sets $\{I_n\}_{n \in \N}$ exhaust all of $\Z$ and since the $\tilde \tau_n$ are bijective and increasing on $I_n$ it follows that the $J_n$ are also a set of subsequent numbers in $\Z$. Thus by injectivity of the $\tilde{\tau}_n$, the $\{J_n\}_{n \in \N}$ exhaust all of $\Z$.
	\end{steps}
\end{Prf}

In a last step we provide a framework, which allows us to pass from a path of metrics to the space of all metrics.

\begin{Def}[discrete family of type (A)]
	\label{DefDiscFamA}
	Let $H$ be a Hilbert space. A discrete family $T:E \to C(H)$ is \emph{self-adjoint of type (A)}, if 
	\begin{enumerate}
		\item There exists a dense subspace $Z \subset H$ such that all $e \in E$ satisfy $\dom T_e = Z$. We set $\dom T:=Z$. 
		\item For all $e \in E$ the operator $T_e$ is self-adjoint.
		\item There exists a norm $| \_ |$ on $Z$ such that for all $e \in E$ the operator $T_e:(Z,| \_ |) \to (H, \| \_ \|_H)$ is bounded and the graph norm of $T_e$ is equivalent to $| \_ |$.
		\item $E$ is a topological space. 
		\item The map $E \to B(Z,H)$, $e \mapsto T_e$, is continuous.
	\end{enumerate}
\end{Def}

\begin{Thm}
	\label{ThmDiscFamA}
	Let $T:E \to C(H)$ be a discrete self-adjoint family of type (A). For any $e_0 \in E$ and any $\varepsilon > 0$ there exists an open neighbourhood $U \subset E$ of $e_0$ such that 
	\begin{align*}
		\forall e \in U: \exists k \in \Z: \forall j \in \Z: d_a(\spc^{e_0}_T(j),\spc^e_T(j+k)) < \varepsilon.
	\end{align*}
\end{Thm}

\begin{Prf}
	Let $\varepsilon > 0$ and $e_0 \in E$ be arbitrary. As in \cref{DefDiscFamA} let $\| \_ \|_H$ be the norm in $H$, $Z:=\dom T$, and let $\| \_ \|_Z$ be the graph norm of $T_{e_0}$ on $Z$. Finally, let $\| \_ \|$ be the associated operator norm in $B(Z,H)$ (which is then also equivalent to the operator norm induced by $| \_ |$).
	
	\begin{steps}
	\stepitem[setup and strategy]
		By construction for any $e_1 \in E$ 
		\begin{align*}
			D_{e_1}(\zeta) := \zeta T_{e_1} + (1-\zeta) T_{e_0} = T_{e_0} + \zeta (T_{e_1} - T_{e_0}), && \zeta \in \C,
		\end{align*}
		defines a discrete self-adjoint holomorphic family $D_{e_1}:\C \to C(H)$ of type (A) with domain $Z$. The idea is to prove the claim using \cref{CorAnalyticGlobalGrowth}. The only problem is that formally for any two $e_1$, $e_2$, the families $D_{e_1}$ and $D_{e_2}$ are different. Hence their constants $C_{I,e_1}$, $C_{I,e_2}$ from \cref{ThmHolAFamilyDerGrowth} for the interval $I$ could differ. Consequently their associated deltas $\delta_{e_1}$, $\delta_{e_2}$ from \cref{CorAnalyticGlobalGrowth} could also differ. We will show that there exists an open neighbourhood $U$ around $e_0$ sufficiently small such that for all $e_1 \in U$ the $\delta = \delta_{e_1}$ is $\geq 1$ , if $\zeta_0$ is always set to $\zeta_0:=t_0:=0$. This will prove the claim.
	\stepitem 
		Recall that by \eqref{EqHolAFamilyDerGrowthDelta} there are $C_1,C_2 > 0$ such that
		\begin{align*} 
			\delta_{e_1} = C_{I,e_1}^{-1} \ln(\min(C_1, \varepsilon C_2)  + 1 ).
		\end{align*} 
		Since $\lim_{t \to 0}{e^t} = 1$, there exists $\varepsilon_1 > 0$ such that 
		\begin{align} \label{PreqDiracSpecContArsinhExpCont}
			\forall |t| \leq 2 \varepsilon_1: \exp(t) - 1 \leq \min(C_1, \varepsilon C_2).
		\end{align}
	\stepitem
		Since $T$ is discrete of type (A), the map $E \to B(Z,H)$, $e \mapsto T_{e}$ is continuous. Consequently, there exists an open neighbourhood $U$ of $e_0$ such that
		\begin{align} \label{EqPreqDiracSpecContArsinhDCont}
			\forall e_1 \in U: \|T_{e_1} - T_{e_0}\| < \min\left( \tfrac{1}{2}, \varepsilon_1 \right).
		\end{align}
	\stepitem
		Now for any $e_1 \in U$, $t \in [0,1]$,  $\varphi \in Z$
		\begin{align*}
			\|D_{e_1}(t) \varphi \|_H
			\geq \|T_{e_0} \varphi \|_H - \| T_{e_1} - T_{e_0} \| \|\varphi \|_Z.
		\end{align*}
		Therefore applying  \eqref{EqHolAFamilyDerGrowthConstantSpec} to $D_{e_1}$,  we obtain 
		\begin{align*}
			\alpha_{I,e_1}
			=\inf_{t \in I}{\inf_{\|\varphi\|_Z=1}}{\|\varphi\|_H + \|D_{e_1}(t) \varphi\|_H} 
			\geq 1 - \| T_{e_0} - T_{e_0} \| 
			\jeq{\eqref{EqPreqDiracSpecContArsinhDCont}}{>} \tfrac{1}{2}.
		\end{align*}
		Furthermore
		\begin{align*}
			\beta_{I,e_1} = \sup_{t \in I}{\|D'_{e_1}(t)\|} = \|T_{e_1} - T_{e_0}\| < \varepsilon_1.
		\end{align*}
		Altogether we achieved for any $e_1 \in U$
		\begin{align*}
			C_{I,e_1} = \alpha_{I,e_1}^{-1} \beta_{I,e_1} < 2 \varepsilon_1.
		\end{align*}
		By \eqref{PreqDiracSpecContArsinhExpCont} this implies
		\begin{align*}
			\exp(C_{I,e_1}) - 1 \leq \min(C_1, \varepsilon C_2)
			\Longrightarrow \delta_{e_1} = C_{I,e_1}^{-1}\ln(\min(C_1, \varepsilon C_2)+1) \geq 1.
		\end{align*}
		This implies the claim.	
	\end{steps}
\end{Prf}

Finally, we apply all our results to Dirac operators.

\begin{Thm}
	\label{ThmSpecContArsinh}
	The map
	\begin{align*}
		\bar \spc:(\Rm(M),\mathcal{C}^1) \to (\Conf,\bar d_a), \; \; g \mapsto \overline{\spc}^g,
	\end{align*}	
	is continous. 
\end{Thm}

\begin{Prf}
	Let $g_0 \in \Rm(M)$ and $\varepsilon > 0$ be arbitrary. By Definition of $\bar d_a$, c.f. \eqref{EqDefbarda}, it suffices to find an open neighbourhood $U \subset \Rm(M)$ such that
	\begin{align} \label{EqSpecContArsinh}
		\forall g' \in U: \exists k \in \Z: \forall j \in \Z: d_a(\spc^{g}(j), \spc^{g'}(j+k)) < \varepsilon.
	\end{align}
	By \cref{ThmSpinorIdentification} the map $\Rm(M) \to B(H^1(\Sigma^{g_0}M),L^2(\Sigma^{g_0}M))$, $h \mapsto \Dirac^h_{g_0}$, is a discrete family of type (A). Consequently, by \cref{ThmDiscFamA} there exists $U$ such that \eqref{EqSpecContArsinh} holds. 
\end{Prf}

In order to apply the Lifting Theorem, we quickly verify that $\pi$ is a covering map.

\begin{Thm}
	\label{ThmPiCoveringArsinh}
	The map $\pi:(\Mon,d_a) \to (\Conf, \bar d_a)$ is a covering map with fibre~$\Z$. 
\end{Thm}

\begin{Prf}
	In this proof we also use \cref{NotTechnical}. By definition of $\tau$, see \eqref{EqDeftau}, $\Z$ acts on $\Mon$ by isometries. In particular $\tau$ is continuous. We will show that for each $u \in \Mon$ there exists an open neighbourhood $V$ such that
	\begin{align} \label{EqTauPropDisc}
		\pi^{-1}(\pi(V)) = \dot{\bigcup}_{j \in \Z}{V.j}.
	\end{align}
	To see this note that the function $\ash \circ u$ is non-decreasing and proper. The set $K_0:=(\ash \circ u)^{-1}(\ash(u(0)))$ is of the form $K_0 = [a_0,b_0]_\Z$ for some $a_0 \leq b_0$, $a_0,b_0 \in \Z$. For the same reason there exist $b_1$, $a_{-1} \in \Z$ such that all in all (see also \Cref{FigCont})
	\begin{align*}
		(\ash \circ u)^{-1}(\ash(u(0))) &= [a_0,b_0]_\Z = K_0, \\
		(\ash \circ u)^{-1}(\ash(u(b_0 + 1))) &= [b_0+1,b_1]_\Z =: K_1, \\
		(\ash \circ u)^{-1}(\ash(u(a_0 - 1))) &= [a_{-1},a_0 - 1]_\Z =:K_{-1}.
	\end{align*}
	Since $\ash(u(\Z))$ is discrete, there exists $\varepsilon > 0$ such that 
	\begin{align*}
		I_{\varepsilon}(\ash(u(0))) \cap I_{\varepsilon}(\ash(u(b_0+1))) = \emptyset, &&
		I_{\varepsilon}(\ash(u(0))) \cap I_{\varepsilon}(\ash(u(a_0-1))) = \emptyset, 
	\end{align*}
	So we obtain open sets
	\begin{align*}
		U_0 := I_{\varepsilon}(\ash(u(0))), && 
		U_1 := I_{\varepsilon}(\ash(u(b_0 + 1))), && 
		U_{-1} := I_{\varepsilon}(\ash(u(a_0-1))), && 
	\end{align*}
	which are mutually disjoint. To see that $V := B_\varepsilon(u)$ satisfies \eqref{EqTauPropDisc} suppose to the contrary that there exists  $v \in V$ and $j \in \Z$ such that $v.j \in V$. Assume $j > 0$ (the proof for $j<0$ is entirely analogous). By hypothesis, this implies that $\ash(v(b_0)) \in U_0$ and $\ash(v(b_0 + j)) = \ash((v.j(b_0))) \in U_0 $. But $\ash \circ v$ is non-decreasing, so $\ash(v(b_0 + j)) \geq \ash(v(b_0 + 1)) \in U_1$. This implies that $\ash(v(b_0+j)) \notin U_0$. Contradiction! \\
	Finally, to see that $\pi$ is a covering map let $u \in [u] \in (\Conf,\bar d_a)$ be arbitrary. Let $V$ be an open neighbourhood of $u$ satisfying \eqref{EqTauPropDisc}. Then $\bar V := \pi(V)$ is evenly covered. Thus $\pi$ is a covering map.

	\begin{figure}[t] 
		\begin{center}
			\includegraphics[scale=0.9, clip=true, trim=90 670 130 70]{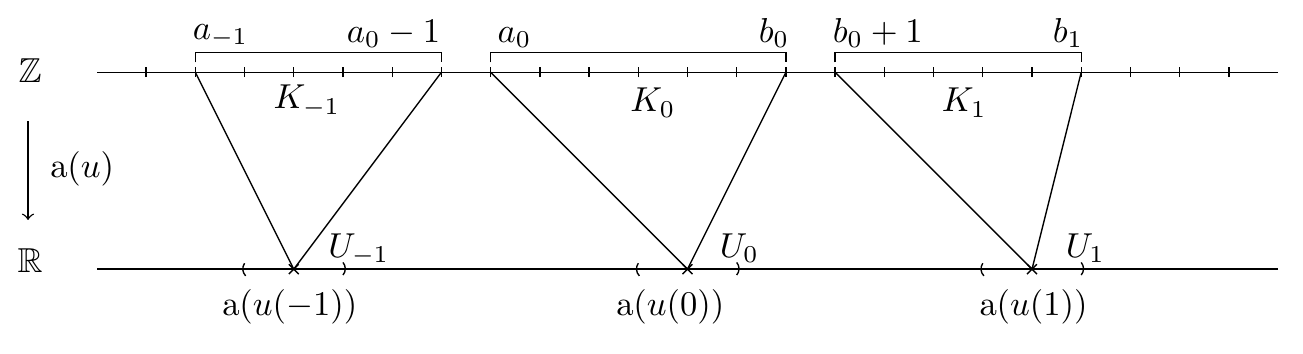}
			\caption{An evenly covered neighbourhood for $u$.}
			\label{FigCont}
		\end{center}
	\end{figure}
\end{Prf}

\section{Moduli Spaces and Spectral Flow}
\label{SecModuliSpecFl}
In this section $M$ is still a compact spin manifold with a fixed topological spin structure $\Theta$ and $I:=[0,1]$ denotes the unit interval. Let $\Diff(M)$ be the diffeomorphism group of $M$. This group acts canonically on the Riemannian metrics via $\Rm(M) \times \Diff(M) \to \Rm(M)$, $(g,f) \mapsto f^*g$. For any subgroup $G \subset \Diff(M)$ the quotient space $\Rm(M) / G$ is called a \emph{moduli space}. We investigate when the map $\bar \spc$ resp. $\widehat{\spc}$ from \cref{MainThmSpec} (and hence the family of functions $(\lambda_j)_{j \in \Z}$ from \cref{MainThmFun}) descends to the moduli spaces, where $G$ is one of the groups 
\begin{align*}
	\Diff^+(M), &&
	\Diff^{\spin}(M), &&
	\Diff^0(M).
\end{align*}
Here $\Diff^+(M)$ denotes the subgroup of orientation-preserving diffeomorphisms and $\Diff^0(M)$ are the diffeomorphisms that are isotopic to the identity. The group $\Diff^{\spin}(M)$ is defined as follows.

\begin{Def}[spin diffeomorphism]
	\label{DefSpinDiffeo}
	A diffeomorphism $f \in \Diff^+(M)$ is a \emph{spin diffeomorphism}, if there exists a morphism $F$ of $\GLtp_m$-fibre bundles  such that
	\begin{align}
		\label{EqDefSpinDiffeo}
		\begin{split}
			\xymatrix{
				\GLtp M \ar[d]^-{\Theta}
					\ar[r]^-{F}
				& \GLtp M
					\ar[d]^-{\Theta}
				\\
				\GLp M
					\ar[d]
					\ar[r]^-{f_*}
				& \GLp M
					\ar[d] 
				\\
				M 
					\ar[r]^-{f}
				& M
			}
		\end{split}
	\end{align}
	commutes. We say $F$ is a \emph{spin lift} of $f$ and define
	\begin{align*}
		\Diff^{\spin}(M)& := \{ f \in \Diff^+(M) \mid \text{$f$ is a spin diffeomorphism} \}.
	\end{align*}
	Notice that if $f$ is a spin diffeomorphism and $M$ is connected, there are precisely two spin lifts $F^\pm$ of $f$ related by $F^+.(-1)=F^-$, where $-1 \in \GLtp_m$.
\end{Def}

\begin{Rem}[spin isometries]
	If $f \in \Diff^{\spin}(M)$, $h \in \Rm(M)$ we can set $g := f^* h$. In this case \eqref{EqDefSpinDiffeo} restricts to the analogous diagram
	\begin{align}	
		\label{EqDefSpinIsometry}
		\begin{split}
			\xymatrix{
				\Spin^g M \ar[d]^-{\Theta^g}
					\ar[r]^-{F}
				& \Spin^h M
					\ar[d]^-{\Theta^h}
				\\
				\SO^g M
					\ar[d]
					\ar[r]^-{f_*}
				& \SO^h M
					\ar[d] 
				\\
				(M,g) 
					\ar[r]^-{f}
				& (M,h)
			}
		\end{split}
	\end{align}
	of metric spin structures. We say that $f$ is a \emph{spin isometry} in this case. Notice that this implies that $(M,g)$ and $(M,h)$ are Dirac-isospectral, i.e. their sets of Dirac eigenvalues are equal as well as the multiplicities of the eigenvalues. Rephrased in the terminology of the previous chapters, this implies $\bar{\spc}^g = \bar{\spc}^h$.
\end{Rem}

These considerations immediately imply the following.
\begin{Thm} 
	\label{ThmPassDiffSpin}
	There exists a commutative diagram
	\begin{align}
		\label{EqDiffSpinPass}
		\begin{split}
			\xymatrix{
				\Rm(M) 
					\ar[d]
					\ar[r]^-{\overline{\spc}} 
				& (\Conf, d_a) 
				\\
				\Rm(M) / \Diff^{\spin}(M)
					\ar@{-->}[ur]^-{\exists \spc^{\spin}}
			}
		\end{split}
	\end{align}	
\end{Thm}

With only a little more work, we get an even stronger statement for $\Diff^0(M)$.

\begin{Thm} 
	\label{ThmPassDiff0}
	There exists a commutative diagram 
	\begin{align}
		\label{EqPassDiff0Mon}
		\begin{split}
			\xymatrix{
				\Rm(M)
					\ar[d]
					\ar[r]^-{\widehat{\spc}} 
				& (\Mon, d_a)
				\\
				\Rm(M) / \Diff^0(M) 
					\ar@{-->}[ur]^-{\exists \spc^0}
			}
		\end{split}
	\end{align}
\end{Thm}

\begin{Prf}  
	The claim will follow from the universal property of the topological quotient, if we can show that 
	\begin{align*}
		\forall g_0 \in \Rm(M): \forall f \in \Diff^0(M): \widehat{\spc}^{g_0} = \widehat{\spc}^{f^*g_0}.
	\end{align*}
	Let $H:M \times I \to M$ be an isotopy from $H_0 = \id$ to $H_1 = f$, let $t \in I$ be arbitrary and set $g_t := H_t^*g_0$. Since $\det((H_t)_*) \neq 0$ and $H_0 = \id$, we obtain $H_t \in \Diff^+(M)$ for all $t \in I$. Consequently, we obtain a diagram
	\begin{align} 
		\label{EqSpinLiftfF}
		\begin{split}
			\xymatrix{
				\GLtp M \times I 
					\ar[d]^-{\Theta}
					\ar@{-->}[r]^-{\exists \tilde H} 
				& \GLtp M \ar[d]^-{\Theta}
				\\
				\GLp M \times I
					\ar[d]
					\ar[r]^-{H_*}
				& \GLp M
					\ar[d] \\
				M \times I
					\ar[r]^-{H} 
				& M.
			}
		\end{split}
	\end{align}
	To be precise, the map $H_*$ is defined by
	\begin{align*}
		\forall B \in \GLp M: \forall t \in I: H_*(B,t) := (H_t)_* B.
	\end{align*}
	To show the existence of $\tilde H$, we notice that since $H$ is an isotopy, it is in particular a homotopy. Consequently $H_* \circ \Theta$ is also a homotopy between $(H_0)_* \circ \Theta = \Theta$ and $(H_1)_* \circ \Theta = f_* \circ \Theta$. Clearly $\widetilde{\id}:\GLtp M \to \GLtp M$ satisfies $\Theta \circ \widetilde{\id} = \Theta = (H_0)_* \circ \Theta$. Since covering spaces have the homotopy lifting property, there exists $\tilde H$ such that $\Theta \circ \tilde H = H_* \circ \Theta$.
	We conclude from \eqref{EqSpinLiftfF} that for any $t \in I$
	\begin{align*}
		\xymatrix{
			\Spin^{g_t} M
				\ar[d]^-{\Theta^{g_t}}
				\ar[r]^-{\tilde H_t}
			& \Spin^{g_0} M 
				\ar[d]^-{\Theta^{g_0}}
			\\
			\SO^{g_t} M
				\ar[d]
				\ar[r]^-{(H_t)_*}
			& \SO^{g_0} M 
				\ar[d]
			\\
			(M,g_t) 
				\ar[r]^-{H_t}
			& (M,g_0)
		}
	\end{align*}
	commutes as well. Consequently, for all $t \in I$, the map $H_t$ is a spin isometry in the sense of \eqref{EqDefSpinIsometry}. Therefore $(M,g_t)$ and $(M,g_0)$ are Dirac isospectral for all $t \in I$. This implies $\widehat{\spc}^{g_0} = \widehat{\spc}^{g_1}$.
\end{Prf}

\begin{Rem}[a counter-example on the torus]
	It remains to discuss the group $\Diff^+(M)$ and one might ask if \eqref{EqDiffSpinPass} still holds, if $\Diff^{\spin}(M)$ is replaced by $\Diff^+(M)$. This is wrong in general. A counter-example can be provided by the standard torus $\T^3 = \R^3 / \Z^3$ with the induced Euclidean metric $\bar g$. It is a well-known fact that the (equivalence classes of) spin structures on $\T^3$ stand in one-to-one correspondence with tuples $\delta \in \Z_2^3$, see for instance \cite{FriedTori}. We denote by $\Spin^{\bar g}_{\delta} \T^3$ the spin structure associated to $\delta$. The map	
	\begin{align*}
		f := 
		\begin{pmatrix}
			1 & 1 & 0\\
			0 & 1 & 0 \\
			0 & 0 & 1 \\
		\end{pmatrix} :\R^3 \to \R^3
	\end{align*}	
	preserves $\Z^3$ and satisfies $\det(f) = 1$. Hence it induces a diffeomorphism $\bar f \in \Diff^+(\T^3)$. One checks that there is a commutative diagram
	\begin{align*}
		\xymatrix{
			\Spin_{(1,1,0)}^{\bar g} \T^3 
				\ar[r]^-{F}
				\ar[d]^{\Theta^{\bar g}} 
			&\Spin^{\bar f^* \bar g}_{(1,0,0)} \T^3
				\ar@{-->}[r]^-{\nexists}
				\ar[d]^{\Theta^{\bar f^* \bar g}}
			& \Spin_{(1,0,0)}^{\bar g} \T^3
				\ar[d]^{\Theta^{\bar g}}
			\\
			\SO^{\bar g} \T^3 
				\ar[r]^-{\bar f^{-1}_*}
				\ar[d]
			&\SO^{\bar f^* \bar g} \T^3
				\ar[r]^-{\bar f_*}
				\ar[d]
			& \SO^{\bar g} \T^3
				\ar[d]
			\\
			(\T^3, \bar g)
				\ar[r]^-{\bar f^{-1}}
			&(\T^3,\bar f^* \bar g)
				\ar[r]^-{\bar f}
			& (\T^3, \bar g) 
		}
	\end{align*}
	The map in the right upper row cannot exist, because otherwise the spin structures corresponding to $(1,1,0)$ and $(1,0,0)$ would be equivalent. The left part of the above diagram is a spin isometry analogous to \eqref{EqDefSpinIsometry}. Therefore $\Dirac^{f^* \bar g}_{(1,0,0)}$ and $\Dirac^{\bar g}_{(1,1,0)}$ are isospectral, but the spectra of $\Dirac^{\bar g}_{(1,1,0)}$ and $\Dirac^{\bar g}_{(1,0,0)}$ are already different as a set. This follows from the explicit computation of the spectra of Euclidean tori, see also \cite{FriedTori}. Consequently $\spec \Dirac^{\bar f^*g}_{(1,0,0)} \neq \spec \Dirac^{\bar g}_{(1,0,0)}$ and no diagram analogous to \eqref{EqDiffSpinPass} can exist for $\Diff^+(\T^3)$.
\end{Rem}

\begin{Rem}
	Notice that in \eqref{EqPassDiff0Mon} the map $\spc^0$ goes from the moduli space for $\Diff^0(M)$ to $\Mon$, where in \eqref{EqDiffSpinPass} the corresponding map $\spc^{\spin}$ goes to $\Conf$. Therefore one might ask, if one could improve \eqref{EqDiffSpinPass} by lifting $\spc^{\spin}$ up to a map $\widehat{\spc}^{\spin}$ such that
\begin{align}
	\label{EqLiftSpinMonDef}
	\begin{split}
		\xymatrix{
			& \Mon
				\ar[d]
			\\
			\Rm(M) / \Diff^{\spin}(M)
				\ar[r]^-{\spc^{\spin}}
				\ar[ur]^-{\widehat{\spc}^{\spin}}
			& \Conf
			}
	\end{split}
\end{align}
	commutes. This question is not so easy to answer and the rest of this section is devoted to the proof that this is not possible in general. To see where the problem is it will be convenient to introduce the following terminology.
\end{Rem}

\begin{Lem}[spectral flow] $ $
    \label{LemSpecFlowDef}
    \begin{enumerate}
        \item 
            For any $f \in \Diff^{\spin}(M)$ and $g \in \Rm(M)$ there exists a unique $\specfl_g(f) \in \Z$ such that
            \begin{align}
                \label{EqDefSpectralFlow}
                \forall j \in \Z: \widehat{\spc}^{g}(j) = \widehat{\spc}^{f^*g}(j - \specfl_g(f)).
            \end{align}
            The induced map $\specfl(f):\Rm(M) \to \Z$ is called the \emph{spectral flow of $f$}. 
        \item
            Let $\mathbf{g}:[0,1] \to \Rm(M)$, $t \mapsto g_t$, be a continuous path of metrics. Let $\spc:\Rm(M) \to \Mon$ be the ordered spectral function for the associated Dirac operators. Take a lift $\hat \spc:\Rm(M) \to \Mon$ of $\overline{\spc}$ such that $\widehat{\spc}^{g_0} = \spc^{g_0}$ as in \eqref{EqMainThmSpc}. There exists a unique integer $\specfl(\mathbf{g}) \in \Z$ such that
            \begin{align*}
                \forall j \in \Z: \widehat{\spc}^{g_1}(j) = \spc^{g_1}(j+\specfl(\mathbf{g})).
            \end{align*}
            The integer $\specfl(\mathbf{g})$ is called the \emph{(Dirac) spectral flow along $\mathbf{g}$}. 
        \item
            \label{ItSpecFlowHomotpyInvariant}
            For any $f \in \Diff^{\spin}(M)$ and any family $\mathbf{g}$ joining $g_0$ and $f^* g_0$, $\specfl_{g_0}(f)=\specfl(\mathbf{g})$.
    \end{enumerate}
\end{Lem}

\begin{Prf} $ $
    \begin{enumerate}
        \item 
            By \cref{ThmPassDiffSpin} the map $\overline{\spc}:\Rm(M) \to \Conf$ descends to a quotient map $\spc^{\spin}:\Rm(M) / \Diff^{\spin}(M) \to \Conf$. This precisely means that $\widehat{\spc}^g$ and $\widehat{\spc}^{f^*g}$ are equal in $\Conf$. By definition of $\Conf$, this implies the existence of $\specfl_g(f)$ as required. 
        \item 
            This follows directly from \eqref{EqMainThmSpc} and the fact that $\spc$ and $\widehat{\spc}$ are equal in $\Conf$.
        \item Set $g_1:=f^*g_0$ and let $\mathbf{g}$ be a family of metrics joining $g_0$ and $g_1$. By \cref{ThmPassDiffSpin}, we obtain $\spc^{g_0} = \spc^{g_1}$. Take a lift $\widehat{\spc}$ satisfying $\widehat{\spc}^{g_0} = \spc^{g_0}$. This implies for all $j \in \Z$
			\begin{align*}
				\spc^{g_0}(j)
				&=\widehat{\spc}^{g_0}(j)
				=\widehat{\spc}^{g_1}(j - \specfl_{g_0}(f)) 
				=\spc^{g_1}(j - \specfl_{g_0}(f) + \specfl(\mathbf{g})) \\
				&=\spc^{g_0}(j - \specfl_{g_0}(f) + \specfl(\mathbf{g})),
			\end{align*}
			which implies $\specfl(\mathbf{g}) - \specfl_{g_0}(f)  = 0$, since $\spc^{g_0}$ is monotonous and all eigenvalues are of finite multiplicity.
    \end{enumerate}
\end{Prf}

\begin{Rem}[spectral flow]
    \label{RemSpecFlowPhillips}
	Intuitively the spectral flow $\specfl(\mathbf{g})$ of a path $\mathbf{g}:[0,1] \to \Rm(M)$ counts the signed number of eigenvalues of the associated path $\Dirac^{g_t}$ of Dirac operators that cross $0$ from below when $t$ runs from $0$ to $1$. The sign is positive, if the crossing is from below, and negative, if it is from above.\\
    The concept of spectral flow is well-known in other contexts. A good introduction can be found in a paper by Phillips, see \cite{phillips}. Phillips introduces the spectral flow for continuous paths $[0,1] \to \mathcal{F}_*^{\sa}$, where $\mathcal{F}_*^{\sa}$ is the non-trivial component of the space of self-adjoint Fredholm operators on a complex separable Hilbert space $H$. In this general setup the definition of spectral flow is a little tricky, see \cite[Prop. 2]{phillips}. But for paths of Dirac operators it coincides with the definition given in \cref{LemSpecFlowDef} above (by \cref{ThmSpinorIdentification} we can think of all Dirac operators $\Dirac^{g_t}$, $t \in [0,1]$, of a path $\mathbf{g}$ as defined on the same Hilbert space). Therefore we have found a convenient alternative to describe the spectral flow in this case using the continuous function $\widehat{\spc}$. \\
    By \cite[Prop. 3]{phillips} the spectral flow of a path of operators depends only on the homotopy class of the path. Consequently, since $\Rm(M)$ is simply-connected, $\specfl(\mathbf{g})$ depends only on $g_0$ and $g_1$. It follows that $\specfl:\Diff^{\spin}(M) \to \Z$ is a group homomorphism.
\end{Rem}

\begin{Rem}
	The map $\widehat{\spc}$ certainly descends to
	\begin{align*}
		\xymatrix{
			\Rm(M)
				\ar[r]^-{\widehat{\spc}}
				\ar[d]
			& \Mon
			\\
			\Rm(M) / \ker \specfl
				\ar@{-->}[ur]^-{\exists}
		}
	\end{align*}
	and $\ker \specfl$ is the largest subgroup of $\Diff^{\spin}(M)$ with this property. Rephrased in these terms we conclude that the map $\spc^{\spin}$ from \eqref{EqDiffSpinPass} lifts to a map $\widehat{\spc}^{\spin}$ as in \eqref{EqLiftSpinMonDef} if and only if $\specfl(f)=0$ for all $f \in \Diff^{\spin}(M)$. Consequently, we have to show the following theorem.
\end{Rem}

\begin{MainThm}
	\label{MainThmFlow}
	There exists a spin manifold $(M,\Theta)$ and a diffeomorphism $f \in \Diff^{\spin}(M)$ such that $\specfl(f) \neq 0$. 
\end{MainThm}

\begin{Prf}
	This proof relies on several other theorems, which are collected in the appendix after this proof for convenient reference (one might want to take a look at these first). The general idea is to obtain $M$ as a fibre $M=P_0$ of a fibre bundle $P \to \S^1$ such that $P$ is spin and $\widehat{A}(P) \neq 0$. This bundle will be isomorphic to a bundle $P_f$ obtained from the trivial bundle $[0,1] \times M \to [0,1]$ by identifying $(1,x)$ with $(0,f(x))$, $x \in M$, for a suitable diffeomorphism $f \in \Diff^{\spin}(M)$. Using various index theorems, we will show that $\widehat{A}(P) = \specfl(f)$. 
	\begin{steps}
	\stepitem[construct a bundle]
		Apply \cref{ThmHanke} to $(k,l) = (1,2)$ and obtain a fibre bundle $P \to S^1$ with some fibre type $M$, where $\dim P = 4n$, $n$ odd, and 
		\begin{align}
			\label{EqIndexHanke}
			\widehat{A}(P) \neq 0.
		\end{align}
		Since $M$ is $2$-connected, $M$ has a unique spin structure. It follows that $m := \dim M = 4n-1 \equiv 3 \mod 4$ and also $m \equiv 3 \mod 8$, since $n$ is odd. Therefore by \cref{ThmAmmannSurgery} there exists a metric $g_0$ on $M$ such that the associated Dirac Operator $\Dirac^M$ is invertible. By \cref{LemClassBundlesS1}, $P$ is isomorphic to $P_f = [0,1] \times M / f$ for some $f \in \Diff^{\spin}(M)$. Define $g_1 := f^* g_0$ and connect $g_0$ with $g_1$ in $\Rm(M)$ by the linear path $g_t := t g_1 + (1-t) g_0$, $t \in [0,1]$. Endow $[0,1] \times M$ with the generalized cylinder metric $dt^2 + g_t$. Denote by $\pi:[0,1] \to \S^1$ the canonical projection. We obtain a commutative diagram
		\begin{align*}
			\xymatrix{
				[0,1] \times M
					\ar[r]
					\ar[d]
				&P_f
					\ar[r]^-{\cong}
					\ar[d]
				&P
					\ar[d]
				\\
				[0,1]
					\ar[r]^-{\pi}
				&\S^1
					\ar[r]^-{\id}
				&\S^1
				}
		\end{align*}
		By construction we can push forward the metric $dt^2 + g_t$ on $[0,1] \times M$ to a metric on $P_f$ and then further to $P$ such that the above row consists of local isometries. The right map is actually an isometry along which we can pull back the spin structure on $P$ to $P_f$. This map is a spin isometry then and therefore we will no longer distinguish between $P$ and $P_f$. The left map is an isometry except that it identifies $\{0\} \times M$ with $\{1\} \times M$. Notice that since $[0,1] \times M$ is simply-connected, the  spin structure on $[0,1] \times M$ obtained by pulling back the spin structure on $P_f$ along $\pi$ is equivalent to the canonical product spin structure on $[0,1] \times M$.

	\stepitem[trivialize]
		The Riemannian manifold $P':= ([0,1] \times M,dt^2 + g_t)$ has two isometric boundary components. Geometrically $P'$ is obtained from $P$ by cutting $M=P_{[0]}$ out of $P$ and adding two boundaries $P'_0$ and $P'_1$, i.e. $P' = (P \setminus P_{[0]}) \coprod P'_0 \coprod P'_1$, where $P'_0$, $P'_1$ are two isometric copies of $P_{[0]}$. By \cref{ThmBB}, we obtain
		\begin{align}
			\label{EqIndexTrivializedBundle}
			\ind(\Dirac_+^P) = \ind(\Dirac_+^{P'}).
		\end{align}
	
	\stepitem[index of half-cylinders]
		Now, set
		\begin{align*}
			\begin{array}{ll}
				Z_0' := (\mathclose] -\infty, 0 \mathclose] \times M, dt^2 + g_0), 
				&Z_1' := (\mathopen[ 1,\infty \mathopen[ \times M, dt^2 + g_1), \\
				Z_1'' := (\mathopen[ 0,\infty \mathopen[ \times M, dt^2 + g_0), 
				&Z := (\R \times M, dt^2+g_0), \\
				Z' := Z_0' \textstyle \coprod Z_1', 
				&Z'' := Z_0' \textstyle \coprod Z_1''.
			\end{array}
		\end{align*}
		Since $Z$ is a Riemannian product, it follows that
		\begin{align}
			\label{EqInvertibleOnProduct}
			\forall \psi \in \Gamma_c(\Sigma Z): \| \Dirac^Z \psi \|_{L^2(\Sigma Z)}^2 \geq \lambda_{\min}^2 \|\psi\|_{L^2(\Sigma Z)}^2,
		\end{align}
		where $\lambda_{\min}$ is the eigenvalue of $\Dirac^{g_0}$ of minimal absolute magnitude. By construction $\Dirac^{g_0}$ is invertible, thus $\lambda_{\min} > 0$. Therefore $\Dirac^{Z}$ is invertible and \emph{coercive at infinity} (see \cref{ThmBB} for the definition). By \cref{ThmBB}, this implies
		\begin{align}
			\begin{split}
				\label{EqIndexHalfZylinderZero}
				0
				&= \ind(\Dirac^{Z}_+)
				= \ind(\Dirac^{Z''}_+)
				= \ind(\Dirac^{Z_0'}_+) + \ind(\Dirac^{Z_1''}_+) \\
				&= \ind(\Dirac^{Z_0'}_+) + \ind(\Dirac^{Z_1'}_+)
				= \ind(\Dirac^{Z'}_+),
			\end{split}
		\end{align}
		where we used the fact that $Z_1'$ and $Z_1''$ are spin isometric.
		
	\stepitem[glue in the half-cylinders]
		Now glue $Z'$ on $P'$ ($Z'_0$ at $\{0\} \times M$ and $Z'_1$ at $\{1\} \times M$) and obtain a bundle $Q =( \R \times M, dt^2 + g_t)$ where $g_t = g_0$ for $t \leq 0$ and $g_t = g_1$ for $t \geq 1$. Since $\Dirac^{g_1}$ is invertible as well, we obtain that  $\Dirac^{Z'}$ satisfies an estimate analogous to \eqref{EqInvertibleOnProduct}. We obtain that $\Dirac^{Q}$ is also coercive at infinity (take $K:=P'$ as the compact subset). By  \cref{ThmBB} we obtain
		\begin{align}
			\label{EqIndexHalfCylinderGlued}
			\ind(\Dirac_+^{Q})
			= \ind(\Dirac_+^{P'}) + \ind(\Dirac_+^{Z'})
			\jeq{\eqref{EqIndexHalfZylinderZero}}{=} \ind(\Dirac_+^{P'}).
		\end{align}
		
	\stepitem[apply hypersurface theory]
		Each $Q_t = \{t\} \times M$ is a hypersurface in $Q$ and $\tilde \partial_t \in \mathcal{T}(Q)$ (horizontal lift of $\partial_t$) provides a unit normal field for all $Q_t$. Therefore we can apply some standard results about the Dirac operator on hypersurfaces, see \cite{BaerGaudMor}: Since $m$ is odd
		\begin{align*}
			\Sigma Q|_{Q_t} = \Sigma^+ Q|_{Q_t} \oplus \Sigma^- Q|_{Q_t} = \Sigma^+ M_t \oplus \Sigma^- M_t, && M_t := (M,g_t),
		\end{align*}
		where $\Sigma^+ M_t = \Sigma^- M_t = \Sigma M_t$ as hermitian vector bundles. The Clifford multiplication $\cdot$ in $\Sigma Q$ is related to the Clifford multiplication $\bullet_t^\pm$ in $\Sigma^\pm M_t$ by $X \bullet_t^\pm \psi = \pm \tilde \partial _t \cdot X \cdot \psi$. Setting
		\begin{align}
			\label{EqDiracCopy}
			\tilde \Dirac^{M_t} := (\Dirac^{M_t} \oplus (-\Dirac^{M_t})),
		\end{align}
		we obtain the Dirac equation on hypersurfaces
		\begin{align}
			\label{EqDiracHyper}
				\tilde \partial_t \cdot \Dirac^{Q} = \tilde \Dirac^{M_t} +\tfrac{m}{2}H_t -  \nabla^{\Sigma Q}_{\tilde \partial_t}: \Gamma(\Sigma Q|_{Q_t}) \to  \Gamma(\Sigma Q|_{Q_t})
		\end{align}
		for all $t \in \R$. 
	
	\stepitem[Identification of the Spinor spaces]
		Let $\psi \in \Gamma(\Sigma Q)$ be a spinor field. For any $t \in \R$ this defines a section a section $\psi_t \in \Gamma(\Sigma Q|_{Q_t})$.  Therefore, we can also think of $\psi$ as a ''section''  in
		\begin{align*}
			\bigcup_{t \in \R}{\Gamma(\Sigma Q|_{Q_t})} \to \R.
		\end{align*}
		and \eqref{EqDiracHyper} tells us how $\tilde \partial_t \cdot \Dirac^Q$ acts on these sections under this identifcation. We would like to apply \cref{ThmSalamon} and therefore have to solve the problem that for various $t$ the Hilbert spaces $\Gamma_{L^2}(\Sigma Q|_{Q_t})$ are different. As discussed in \cite{BaerGaudMor}, we will use the following identification: For any $x \in M$, consider the curve $\gamma^x:\R \to \R \times M$, $t \mapsto (t,x)$. Each spinor $\psi \in \Gamma(\Sigma Q)$ determines a section $\psi^x$ in $\Sigma Q$ along $\gamma_x$. Using the connection $\nabla^{\Sigma Q}$, we obtain a parallel translation $\tau^{t}_{0}:\Sigma_x M_{t} \to \Sigma_x M_{0}$. Notice that the Clifford multiplication ''$\cdot$'', the vector field $\tilde \partial_t$ and the volume form $\omega$ that determines the splitting $\Sigma Q = \Sigma^+Q \oplus \Sigma^-Q$ are all parallel.
		Identifying $M$ with $M_0$ we obtain a map
			\DefMap{\tau(\psi):=\bar \psi: \R}{\Gamma(\Sigma^+ M) \oplus \Gamma(\Sigma^- M)}{t}{(x \mapsto \tau_{t}^0(\psi^+_{(t,x)})+\tau_{t}^0(\psi^-_{(t,x)}))}
		This identification defines an isometry $\tau: \Gamma_{L^2}(\Sigma Q) \to L^2(\R, \Gamma_{L^2}(\Sigma^+ M \oplus \Sigma^- M))$. The operator in \eqref{EqDiracHyper} can be pulled back via a commutative diagram
		\begin{align*}
			\xymatrix{
				\Gamma_{L^2}(\Sigma Q)
					\ar[r]
					\ar[d]^-{\tau}
				&\Gamma_{L^2}(\Sigma Q)
					\ar[d]^-{\tau}
				\\
				L^2(\R, \Gamma_{L^2}(\Sigma^+ M \oplus \Sigma^- M))
					\ar[r]
				&L^2(\R, \Gamma_{L^2}(\Sigma^+ M \oplus \Sigma^- M))
				}
		\end{align*}
		and the upper row is given by \eqref{EqDiracHyper}. We calculate how this equation looks like in the lower row: Since Clifford multiplication and $\tilde \partial_t$ are parallel
		\begin{align*}
			 \tau \circ (\tilde \partial_t \cdot \Dirac^{Q}) \circ \tau^{-1}
			 = \tilde \partial_t \cdot (\tau \circ \Dirac^{Q} \circ \tau^{-1})
			 =: \tilde \partial_t \cdot \Dirac^{Q}_{M}.
		\end{align*}
		Since there is a splitting $(\Sigma Q, \nabla^Q) = (\Sigma^+Q,\nabla^{+} \oplus \nabla^{-})$  it suffices to check the following for a $\psi \in \Gamma(\Sigma^+ Q)$: Let $x \in M$ be arbitrary and let $D_x$ be the covariant derivative induced by $\nabla^{\Sigma Q}$ along $\gamma^x$. For any $t_0 \in \R$ 
		\begin{align*}
			\nabla_{\tilde \partial t}^{\Sigma Q}{\psi}|_{{(t_0,x)}}
			&=D_x(\psi^x)(t_0)
			=\lim_{t \to t_0}{\frac{\tau_t^{t_0}(\psi^x(t)) - \psi^x(t_0)}{t-t_0}}
		\end{align*}
		and consequently
		\begin{align*}
			\overline{\nabla_{\tilde \partial t}^{\Sigma Q}{\psi}}(t_0)|_x
			=&\tau_{t_0}^0(\nabla_{\tilde \partial t}^{\Sigma Q}{\psi}|_{{(t_0,x)}})
			=\tau_{t_0}^0\left( \lim_{t \to t_0}{\frac{\tau_t^{t_0}(\psi^x(t)) - \psi^x(t_0)}{t-t_0}} \right) \\
			=& \lim_{t \to t_0}{\frac{\tau_{t_0}^0(\tau_t^{t_0}(\psi^x(t))) - \tau_{t_0}^0(\psi^x(t_0))}{t-t_0}} 
			= \lim_{t \to t_0}{\frac{\tau_t^{0}(\psi_{(t,x)}) - \tau_{t_0}^0(\psi|_{(t_0,x)})}{t-t_0}} \\
			=& \lim_{t \to t_0}{\frac{\bar \psi(t)|_x - \bar \psi(t_0)|_x}{t-t_0}} 
			=\tfrac{d\bar \psi}{dt}(t_0)|_x.
		\end{align*}
		This implies for any $\bar \psi \in L^2(\R, \Gamma_{H^1}(\Sigma^+ M \oplus \Sigma^- M))$
		\begin{align*}
			(\tau \circ \nabla_{\tilde \partial t}^{\Sigma Q} \circ \tau^{-1})(\bar \psi) = \tfrac{d}{dt} \bar \psi.
		\end{align*}		
		All in all \eqref{EqDiracHyper} transforms under $\tau \circ \_ \circ \tau^{-1}$ into
	\begin{align}
			\label{EqDiracHyperTransform}
				\tilde \partial_t \cdot \Dirac^{Q}_M = \tilde \Dirac^{Q}_M + \tfrac{m}{2}H - \tfrac{d}{dt},
		\end{align}
		where $H:\R \to \R$, $t \mapsto H_t$ and $\tilde \Dirac^Q_M:\R \to \Gamma(\Sigma M \oplus \Sigma M)$ is given by $\tilde \Dirac^Q_M(t) = \bar \Dirac(t) \oplus (- \bar \Dirac(t))$, $\bar \Dirac(t) = \tau \circ \Dirac^{M_t} \circ \tau^{-1}$. 
	
	\stepitem[apply Salamon] Set $H := L^2(\R, \Gamma_{L^2}(\Sigma M))$, $W := L^2(\R, \Gamma_{H^1}(\Sigma M))$ and $A(t) = \bar \Dirac(t)$. Using \cref{ThmSalamon}, we obtain 
	\begin{align}
		\label{EqIndexSpecFlowSalamon}
		\begin{split}
			\ind(\Dirac^Q_+)
			&=\ind((\tilde \partial_t \cdot \Dirac^Q_M - \tfrac{m}{2} H)_+) 
			\jeq{\eqref{EqDiracHyperTransform}}{=}\ind((\Dirac^{Q}_M  - \tfrac{d}{dt})_+) \\
			&=\ind((\bar \Dirac  - \tfrac{d}{dt})_+) 
			=\ind(\tfrac{d}{dt} - \bar \Dirac)
			=\specfl(\bar \Dirac).
		\end{split}
	\end{align}
	Now the spectral flow $\specfl(\bar \Dirac)$ in the sense of Salamon, c.f. \cite{Salamon}, coincides with the spectral flow in the sense of \cref{LemSpecFlowDef}.
	
	\stepitem[final argument]
		By the classical Atiyah-Singer Index theorem, c.f. \cite[Thm. III.13.10]{LM}, we obtain 
		\begin{align}
			\label{EqIndexAtiyahSingerClassical}
			\widehat{A}(P) = \ind(\Dirac_+^P). 	
		\end{align}
		Consequently, we can put all the steps together to obtain
		\begin{align*}
			0 \jeq{\eqref{EqIndexHanke}}{\neq} \widehat{A}(P)
			&\jeq{\eqref{EqIndexAtiyahSingerClassical}}{=} \ind(\Dirac_+^P)
			\jeq{\eqref{EqIndexTrivializedBundle}}{=} \ind(\Dirac_+^{P'}) \\
			&\jeq{\eqref{EqIndexHalfCylinderGlued}}{=} \ind(\Dirac_+^{Q}) 
			\jeq{\eqref{EqIndexSpecFlowSalamon}}{=} \specfl(\bar \Dirac)
			=\specfl_{g_0}(f).
		\end{align*}
	\end{steps}
\end{Prf}

\begin{appendix}

\section{Appendix}
Here we collect some results for convenient reference.

\begin{Thm}[{\cite[Thm. 1.3]{Hanke}}]
	\label{ThmHanke}
	Given $k,l \geq 0$ there is an $N=N(k,l) \in \N_{\geq 0}$ with the following property: For all $n \geq N$, there is a $4n$-dimensional smooth closed spin manifold $P$ with non-vanishing $\widehat{A}$-genus and which fits into a smooth fibre bundle
	\begin{align*}
		X \to P \to \S^k.
	\end{align*}
	In addition we can assume that the following conditions are satisfied:
	\begin{enumerate}
		\item The fibre $X$ is $l$-connected.
		\item The bundle $P \to \S^k$ has a smooth section $s:\S^k \to P$ with trivial normal bundle.
	\end{enumerate}
\end{Thm}

\begin{Thm}[{\cite[Thm 1.1]{AmmannSurgery}}]
	\label{ThmAmmannSurgery}
	Let $M$ be a compact spin manifold of dimension $m$. Assume $m \not \equiv 0 \mod 4$ and $m \not \equiv 1,2 \mod 8$. Then there exists a metric $g \in \Rm(M)$ such that $\Dirac^g$ is invertible. 
\end{Thm}

\begin{Lem}[Classification of fibre bundles over $\S^1$]
	\label{LemClassBundlesS1}
	Let $M$ be any smooth manifold.
	\begin{enumerate}
		\item For any $f \in \Diff(M)$ the space obtained by setting
		\begin{align*}
			P_f := [0,1] \times M / \sim, &&
			\forall (t,x) \in [0,1] \times M: (1,x) \sim (0,f(x))
		\end{align*}
		is a smooth $M$-bundle over $\S^1 = [0,1] / (0 \sim 1)$. 
		\item Let $\mathcal{I}$ denote the isotopy classes of $\Diff(M)$ and $\mathcal{P}$ the isomorphism classes of $M$-bundles over $\S^1$. The map
			\DefMap{\mathcal{I}}{\mathcal{P}}{{[f]}}{{[P_f]}}
		is well-defined, surjective and if $[P_{f}] = [P_{f'}]$, then $f$ is isotopic to a conjugate of~$f'$. 

		\item Let $M$ be oriented. Then $P_f$ is orientable if and only if $f \in \Diff^+(M)$.
		\item Let $M$ be spin and simply-connected. Then $P_f$ is spin if and only if $f \in \Diff^{\spin}(M)$. 
	\end{enumerate}
\end{Lem} 

The proof of this is elementary. 

\begin{Thm}[{\cite[Thm. 8.17, Rem. 8.18a)]{BaerBallmann}}] 
	\label{ThmBB}
	Let $(M,g)$ be an even-dimensional smooth complete spin manifold with volume element $\mu$ and a fixed spin structure. Let $N$ be a closed and two-sided hypersurface in $M$. Cut $M$ open along $N$ to obtain a manifold $M'$ with two isometric boundary components $N_1$ and $N_2$. Consider the pullbacks $\mu'$, $\Sigma M'$, $\Dirac'$ of $\mu$, $\Sigma M$ and $\Dirac$. Let $\Dirac$ be \emph{coercive at infinity}, i.e. there exists a compact subset $K \subset M$ and a $C>0$ such that 
	\begin{align*}
		\|\psi\|_{L^2(\Sigma M)} \leq C \| \Dirac \psi \|_{L^2(\Sigma M)}
	\end{align*}
	for all $\psi \in \Gamma(\Sigma M)$, which are compactly supported in $M \setminus K$. Then $\Dirac'$ is Fredholm and
	\begin{align*}
		\ind \Dirac'_+ = \ind \Dirac_+.
	\end{align*}
	Here $\Dirac'$ is to be understood as the Dirac operator with APS-boundary conditions, i.e.
	\begin{align*}
		\dom(\Dirac') = \{ \psi \in H^1(\Sigma M ') \mid \psi|_{N_1} \in H^{1/2}_{\mathclose]-\infty,0\mathclose]}(\tilde \Dirac), \psi|_{N_2} \in H^{1/2}_{\mathopen[0, \infty \mathopen[}(\tilde \Dirac)\},
	\end{align*}
	where $\tilde \Dirac $ is the Dirac operator on the boundary (resp. its two-fold copy as in \eqref{EqDiracCopy}).
\end{Thm}

\begin{Thm}[{\cite[Thm A]{Salamon}}]
	\label{ThmSalamon}
	Assume we are given the following data.
	\begin{enumerate}
		\item $(H,\|\_\|_H)$ be a complex separable Hilbert space.
		\item Let $W \subset H$ be a dense subspace, $\| \_ \|_W$ be a norm on $W$ such that $W$ is also a Hilbert space and the injection $W \hookrightarrow H$ is compact.
		\item $\{A(t)_{t\in \R}\}$ is a family of unbounded self-adjoint operators on $H$ with time independent domain $W$ such that for each $t \in \R$ the graph norm of $A(t)$ is equivalent to $\| \_ \|_W$.
		\item The map $\R \to L(W,H)$, $t \mapsto A(t)$, is continuously differentiable with respect to the weak operator topology.
		\item There exist invertible operators $A^\pm \in L(W,H)$ such that $\lim_{t \to \pm \infty}{A(t)} = A^\pm$ in norm topology.
	\end{enumerate}
	 Then the operator
	 \begin{align*}
		D_A := \tfrac{d}{dt} - A(t): W^{1,2}(\R,H) \cap L^2(\R,W) \to L^2(\R,H)
	 \end{align*}
	 is Fredholm and its Fredholm index is equal to the spectral flow $\specfl(A)$ of the operator family $A=(A(t))_{t \in \R}$. 
\end{Thm}

\end{appendix}

\vspace{3em}

\textbf{Acknowledgements.} I would like to thank my PhD supervisor Bernd Ammann very much for his continuing support. I am also grateful to our various colleagues at the University of Regensburg, in particular Ulrich Bunke, Nicolas Ginoux and Andreas Hermann for fruitful discussions. Furthermore I am indebted to Nadine Grosse for explaining \cite{BaerBallmann} to me. This research was enabled by the \emph{Studienstiftung des deutschen Volkes}.

\addcontentsline{toc}{section}{References}

\printbibliography

\end{document}